\documentclass[smallcondensed]{svjour3}     
\smartqed  
\usepackage{threeparttable}
\usepackage{enumerate}
\usepackage{algorithm}
\usepackage{algorithmic}
\usepackage{graphicx}
 \usepackage[misc]{ifsym}
\usepackage{graphicx,amsmath,amssymb,txfonts}
\usepackage{subfigure}
\usepackage{epsfig}
\usepackage{hypernat}
\usepackage{multirow}
\usepackage{booktabs}
\usepackage[colorlinks=true,linkcolor=blue,citecolor=blue]{hyperref}

\begin{document}

\title{An extrapolation cascadic multigrid method combined with a fourth-order compact scheme for 3D Poisson equation}

\titlerunning{An EXCMG method combined with a fourth-order compact scheme for 3D Poisson equation}        

\author{ Kejia Pan\and Dongdong He      \and Hongling Hu
}


\institute{K.J. Pan   \at
              School of Mathematics and Statistics, Central South University, Changsha 410083, China\\
              \email{pankejia@hotmail.com}
           \and
            D.D. He (\Letter) \at
              School of Aerospace Engineering and Applied Mechanics,
Tongji University, Shanghai 200092, China\\
              \email{dongdonghe@tongji.edu.cn}
          \and
          H.L. Hu \at
          College of Mathematics and Computer Science, Key Laboratory of High Performance Computing and Stochastic Information Processing (Ministry of Education of China), Hunan Normal University, Changsha 410081, China   \\       
           \email{hhling625@163.com}
}

\date{Received: date / Accepted: date}

\maketitle

\begin{abstract}
Extrapolation cascadic multigrid (EXCMG) method is an efficient multigrid method which has mainly been used for solving the two-dimensional elliptic boundary value problems with linear finite element discretization in the existing literature. In this paper, we develop an EXCMG method to solve the three-dimensional Poisson equation  on rectangular domains by using the compact finite difference (FD) method with unequal meshsizes in different coordinate directions. The resulting linear system from compact FD discretization is solved by the conjugate gradient (CG) method with a relative residual stopping criterion. By combining the Richardson extrapolation and tri-quartic Lagrange interpolation for the numerical solutions from two-level of grids (current and previous grids), we are able to produce an extremely accurate approximation of the actual numerical solution on the next finer grid, which can greatly reduce the number of relaxation sweeps needed. Additionally, a simple method based on the midpoint extrapolation formula is used for the fourth-order FD solutions on two-level of grids to achieve sixth-order accuracy on the entire fine grid cheaply and directly. The gradient of the numerical solution can also be easily obtained through solving a series of tridiagonal linear systems resulting from the fourth-order compact FD discretizations. Numerical results show that our EXCMG method is much more efficient than the classical V-cycle and W-cycle multigrid methods. Moreover, only few CG iterations are required on the finest grid to achieve full fourth-order accuracy in both the $L^2$-norm and $L^{\infty}$-norm for the solution and its gradient when the exact solution belongs to $C^6$. Finally, numerical result shows that our EXCMG method is still effective when the exact solution has a lower regularity, which widens the scope of applicability of our EXCMG method.

\end{abstract}

\section{Introduction}
Poisson equation is a partial differential equation of elliptic type with broad application in electrostatics, mechanical engineering, theoretical physics and geophysics. The Dirichlet boundary value problem for the three-dimensional (3D) Poisson equation has the following form:
\begin{equation}\label{bvp}
\left\{ \begin{aligned}
         u_{xx}+u_{yy}+u_{zz}&=&f(x,y,z),\quad  &\textrm{in } \Omega,\\
u(x,y,z)&=&g(x,y,z),\quad &\textrm{on } \partial\Omega,
        \end{aligned} \right.
\end{equation}
where $\Omega$ is a 3D rectangle domain and $\partial\Omega$ is its boundary.
Here we assume that the forcing function $f(x,y,z)$, the boundary function $g(x,y,z)$ and the exact solution $u(x,y,z)$
are continuously differentiable and have the necessary continuous partial derivatives up to certain orders.

 The compact finite difference (FD) method for solving Poisson equations has been well studied since 1984~\cite{Strikwerda,Gupta1,Gupta2,Spotz,Sutmann,Wang1,Gupta3,Othman,Schaffer,Zhang1,Zhang2,Zhang3,Ge}.
Specifically, two-dimensional (2D) and 3D  Poisson equations can be solved by high-order compact FD methods~\cite{Strikwerda,Gupta1,Gupta2,Spotz,Sutmann,Wang1}.
These schemes are called ``compact'' since they only use minimum  grid points to achieve fourth-order accuracy explicitly in the discretization formulas. Moreover, there has been a renewed interest in combining high-order compact scheme with multigrid method to solve Poisson equations. The classical multigrid method~\cite{McCormick,Briggs,Trottenberg} combined with compact FD method for solving 2D and 3D Poisson equations has been conducted in \cite{Gupta3,Othman,Schaffer,Zhang1,Zhang2,Zhang3,Ge,Moghaderi}. For example, Wang and Zhang~\cite{Zhang2} proposed a Richardson extrapolation for the numerical solutions from the two-level grids together with an operator based interpolation iterative strategy to achieve sixth-order accuracy by using the classical multigrid method and the fourth-order compact FD scheme.  Ge~\cite{Ge} developed a fourth-order compact FD method with the classical multigrid method to solve the 3D Poisson equation using unequal meshsizes in different coordinate directions.  Dehghan et al.~\cite{Moghaderi} solved the 1D, 2D and 3D Poisson equations with both second-order and fourth-order compact FD methods based on a new two-grid multigrid method. Besides Poisson equation, the classical multigrid method  has been applied to many problems, including the biharmonic equation~\cite{Altas}, the convection-diffusion equation~\cite{Zhang33,Ge2,Wang2010} and so on.

Cascadic multigrid (CMG) method  proposed by Deuflhard and Bornemann in~\cite{Bornemann} is a variant of the multigrid without any coarse grid correction steps, where
instead of starting from the finest grid, the solution is first computed on the coarsest grid and the recursively interpolated and
relaxed on finer grids. Bornemann and Deuflhard \cite{Bornemann} showed that it is an optimal iteration method with respect to the energy norm.
Since the 1990s, the method has been frequently used to solve the elliptic equation with the finite element (FE) discretization
because of its high efficiency and simplicity~\cite{Braess1, Timmermann,Shaidurov2, Shaidurov3, Shi1, Braess2, Stevenson, Zhou,  Du,  Xu, Yu}.
In 2007,  Shi et al.~\cite{Shi07} proposed an economical cascadic multigrid method using the different criteria for choosing the smoothing steps on each level of grid. Later,  based on a new Richardson extrapolation formula for the linear FE solution, an extrapolation cascadic multigrid (EXCMG) method was first proposed by Chen et al.~\cite{Chen1,Chen2} to solve 2D Poisson equation with the linear FE discretization.  For the EXCMG method, in order to obtain a better initial guess of the iterative solution on the next finer grid, numerical solutions on the two-level of grids (current and previous grids) are needed (whereas only one-level of numerical solution is needed in the CMG method). The EXCMG algorithm has been successfully applied to non-smooth problems~\cite{Hu2014}, linear parabolic problems~\cite{Hu2014b}, and the simulation of the electric field with a point singularity arising in geophysical exploration~\cite{Pan2012,Pan2014}.
However, to our best knowledge,  the EXCMG algorithm has mainly been used for solving the 2D elliptic problems with the linear FE discretization in existing literature. But it is of more importance to solve the 3D problems efficiently and accurately arising in many engineering areas, such as geophysical exploration~\cite{Newman2014}. Since the construction process of the higher-order (at least fifth-order) approximation to the fourth-order compact FD solution on the next finer grid has to be different from  the construction process of the third-order approximation to the second-order FE solution, it will be nontrivial to extend the EXCMG method from 2D problems with second-order FE discretization to 3D problems with fourth-order compact FD discretization.

In this paper, we will propose an EXCMG method combined with the fourth-order compact difference scheme  to solve the Dirichlet boundary value problem of the 3D Poisson equation (\ref{bvp}) in rectangular domains. In our approach, the computational domain is discretized by  regular grids, and a 19-point fourth-order compact difference scheme  is used to discretize the 3D Poisson equation with unequal meshsizes in different directions.
By combining the Richardson extrapolation and tri-quartic Lagrange interpolation for the numerical solutions from two-level of grids (current and previous grids), we are able to obtain a much better initial guess of the iterative solution on the next finer grid than one obtained by using linear interpolation in CMG method. Then, the resulting large linear system is solved by the conjugate gradient (CG) solver using the above obtained initial guess. Additionally, a tolerance related to relative residual is introduced in the CG solver in order to obtain conveniently the numerical solution with the desired accuracy. Moreover, when the exact solution is sufficiently smooth, a simple method based on the midpoint extrapolation formula can be  used to obtain cheaply and directly a sixth-order accurate solution on the entire fine grid from two fourth-order FD solutions on two different scale grids (current and previous grids). And a fourth-order compact FD scheme can be used to compute the gradient of the solution by solving a series of tridiagonal linear systems. Finally, our method has been used to solve 3D Poisson equations with more than 16 million unknowns in about 10 seconds on a desktop with 16GB RAM installed, which is much more efficient than the classical multigrid methods. 

The rest of the paper is organized as follows:
section~\ref{CFD} gives the description of the compact FD discretization for the 3D Poisson equation. Section~\ref{mg} reviews the classical V-cycle and W-cycle multigrid methods. In section \ref{ecmg}, we first derive some sixth-order extrapolation formulas, and then develop a new EXCMG method to solve 3D Poisson equation. Section~\ref{sec5} presents the numerical results to demonstrate the high efficiency and accuracy of the proposed method. And conclusions are given in the final section.

\section{Compact difference scheme}~\label{CFD}

We consider a cubic domain $\Omega=[0,L_x]\times[0, L_y]\times[0, L_z]$, and discretize the domain with unequal meshsizes $h_x, h_y$ and $h_z$ in the $x, y$ and $z$
 coordinate directions, respectively. Let $N_x = L_x/h_x$, $N_y = L_y/h_y$, $N_z = L_z/h_z$ be the numbers of uniform
intervals along the $x$, $y$ and $z$ directions. The grid points are ($x_i,y_j, z_k$), with $x_i = ih_x, y_j = jh_y$ and $z_k = kh_z, i = 0,1,\cdots ,N_x, j = 0,1,\cdots,N_y$ and $k = 0,1,\cdots,N_z$. The quantity $u_{i,j,k}$ represents the numerical solution at ($x_i,y_j, z_k$).

Then the value on the boundary points $u_{i,j,k}(i=0,N_x \textrm{ or } j=0,N_y \textrm{ or } k=0,N_z)$ can be evaluated directly  from the Dirichlet boundary condition.
For internal grid points ($i=1,\cdots,N_x-1, j=1,\cdots,N_y-1, k=1,\cdots,N_z-1$), the 19-point  fourth-order compact difference scheme with unequal-meshsize for 3D Poisson equation was derived in~\cite{Wang1,Ge}:
\begin{align}\label{compact}
&-8\left(\frac{1}{h^2_x}+\frac{1}{h^2_y}+\frac{1}{h^2_z}\right)u_{i,j,k}+(\frac{4}{h^2_x}-\frac{1}{h^2_y}-\frac{1}{h^2_z})\left(u_{i+1,j,k}+u_{i-1,j,k}\right)+(\frac{4}{h^2_y}-\frac{1}{h^2_x}-\frac{1}{h^2_z})\left(u_{i,j+1,k}+u_{i,j-1,k}\right)\nonumber\\
&+\left(\frac{4}{h^2_z}-\frac{1}{h^2_x}-\frac{1}{h^2_y}\right)\left(u_{i,j,k+1}+u_{i,j,k-1}\right)+\frac{1}{2}\left(\frac{1}{h^2_x}+\frac{1}{h^2_y}\right)\left(u_{i+1,j+1,k}+u_{i+1,j-1,k}+u_{i-1,j+1,k}+u_{i-1,j-1,k}\right)\nonumber\\
&+\frac{1}{2}\left(\frac{1}{h^2_x}+\frac{1}{h^2_z}\right)\left(u_{i+1,j,k+1}+u_{i+1,j,k-1}+u_{i-1,j,k+1}+u_{i-1,j,k-1}\right)\nonumber\\
&+\frac{1}{2}\left(\frac{1}{h^2_y}+\frac{1}{h^2_z}\right)\left(u_{i,j+1,k+1}+u_{i,j-1,k+1}+u_{i,j+1,k-1}+u_{i,j-1,k-1}\right)\nonumber\\
&=\frac{1}{2}(6f_{i,j,k}+f_{i+1,j,k}+f_{i-1,j,k}+f_{i,j+1,k}+f_{i,j-1,k}+f_{i,j,k-1}+f_{i,j,k+1}).
\end{align}

Let $h=\max\{h_x, h_y, h_z\}$, throughout this paper, we denote $u_h$ to be the FD solution of (\ref{compact}) with mesh sizes $h_x, h_y, h_z$, while use $u_{h/2}$ to denote the FD solution  of (\ref{compact}) when mesh sizes are $h_x/2, h_y/2, h_z/2$.   Then the difference scheme (\ref{compact}) can be expressed in the following matrix form:
\begin{equation}\label{FD}
  A_h u_h = f_h,
\end{equation}
where $A_h$ is  a sparse positive definite matrix, and  $f_h$ denotes the right hand-side vector of (\ref{compact}) with mesh sizes $h_x, h_y$ and $h_z$. 

\section{Classical multigrid methods}\label{mg}

\begin{figure}
  \centering
  \includegraphics[width=4in]{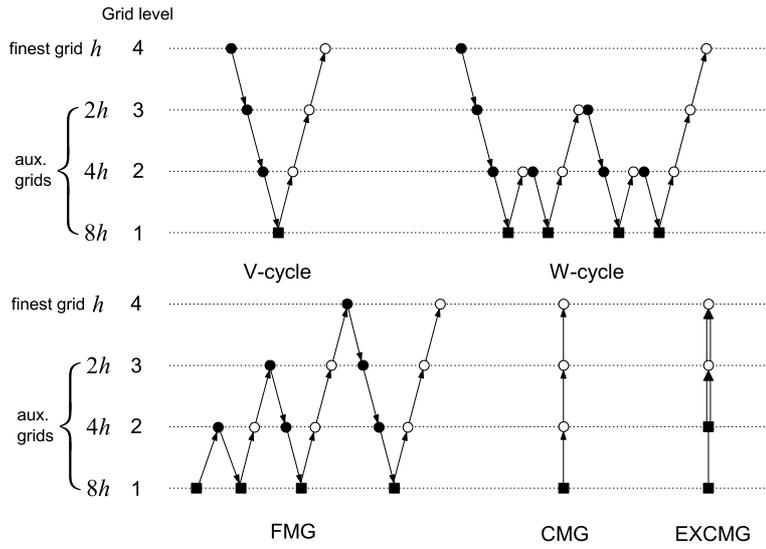}\\
  \caption{The four-level structure of the V-cycle, W-cycle, FMG, CMG and EXCMG methods. In the diagram, $\bullet$ denotes pre-smoothing,
  $\circ$ denotes post-smoothing, $\uparrow$ denotes prolongation (usually defined by linear interpolation), $\downarrow$ denotes restriction,
  $\Uparrow$ denotes extrapolation and high-order interpolation, and $\blacksquare$ denotes direct solver.
  }\label{VW}\end{figure}

The multigrid method is based on the idea that classical relaxation methods
strongly damp the oscillatory error components, but converge slowly for
smooth error components~\cite{Briggs,Trottenberg}.
Hence, after a few relaxation sweeps, we compute
the smooth residual of the current approximation $v_h$ (with mesh sizes $h_x,h_y,h_z$) and transfer it to a coarser grid $\Omega_{2h}$ (with mesh sizes $2h_x,2h_y,2h_z$) by a restriction operation, where the errors become more oscillatory.
Solving the residual equation on the
coarse grid $\Omega_{2h}$, interpolating the correction back to the fine grid $\Omega_{h}$, and adding
it to the fine-grid current approximation $v_h$ yields to the two-grid correction method.
Since the coarse-grid problem is not much different from the original problem,
we can perform a few, say $\gamma$,  two-grid iteration steps (see Fig. \ref{VW}) to the residual equation on the coarse grid, which
means relaxing there and then moving to $\Omega_{4h}$ (with mesh sizes $4h_x,4h_y,4h_z$) for the correction step. We can repeat
this process on successively coarser grids until a direct solution of the residual
equation is possible. Then the corrections are interpolated back
to finer grids until the process reaches the finest grid $\Omega_h$ (with mesh sizes $h_x,h_y,h_z$) and the fine-grid
approximate solution is corrected.

Usually, the cases $\gamma= 1$ and $\gamma = 2$ are particularly interesting. We refer to the case $\gamma = 1$ as V-cycle and to $\gamma = 2$ as
W-cycle. The number $\gamma$ is also called cycle index.
A V-cycle multigrid method is obtained when the V-cycle is repeated until a stopping criterion is satisfied on the finest grid.
We refer to a V-cycle (W-cycle) with $\nu_1$ relaxation sweeps before the correction step and $\nu_2$ relaxation
sweeps after the correction step as a V($\nu_1$, $\nu_2$)-cycle (W($\nu_1$, $\nu_2$)-cycle).

\section{Extrapolation cascadic multigrid methods}\label{ecmg}

The CMG method proposed by Deuflhard and Bornemann in \cite{Bornemann} is a variant of full multigrid (FMG) method without any coarse grid correction steps but with an a posteriori control of the number of smoothing iterations (see Fig. \ref{VW}). It has been shown that the CMG method has optimal computational complexity for both conforming and nonconforming elements with CG as a smoother.
Since the 1990s, the CMG method has received quite a bit of attention from researchers because of its high efficiency and simplicity~\cite{Shaidurov1,Braess1, Timmermann,Shaidurov2, Shaidurov3, Shi1, Braess2, Stevenson, Zhou,  Du,  Xu, Yu}.

In 2008, by using Richardson extrapolation and bilinear quadratic interpolation for the FE solutions on two-level of grids (current and previous grids) to obtain an extremely accurate initial guess of the iterative solution on the next finer grid, Chen et al.~\cite{Chen1} proposed an extrapolation cascadic multigrid (EXCMG) method to solve 2D elliptic boundary value problems. It has been shown in~\cite{Chen2} that the EXCMG method is much more efficient than the CMG method, which simply uses the linear interpolation for the FE solution on the current grid to provide an initial guess of the iterative solution on the next finer grid.
Recently, we improved and generalized the EXCMG method to solve large linear systems resulting from FE discretization of 3D elliptic problems,
compared it with the classical multigrid methods, and further presented the reason why EXCMG algorithms are highly efficient~\cite{Pan2015}. However, to our best knowledge, CMG and EXCMG are mainly  used for linear FE method in existing literature, and it will be interesting to extend the EXCMG method to the field of high-order FD method.

\subsection{EXCMG algorithm combined with compact difference scheme}\label{subsec}

The key ingredients of the EXCMG method are extrapolation and high-order interpolation (see Fig. \ref{VW}), which can produce a much better initial guess of the iterative solution on the next finer grid than one obtained by using linear interpolation in CMG method.

In this subsection, we will propose a new EXCMG method combined with fourth-order compact difference scheme for solving the Dirichlet boundary value problem of the 3D Poisson equation, which is stated in the following algorithm.

\begin{algorithm}[!htb]         

\caption{New EXCMG method: $(u_h, \tilde{u}_h)$ $\Leftarrow$ EXCMG($A_h, f_h, L ,\epsilon$)}             

\label{alg:oEXCMG}                  

\begin{algorithmic}[1]                

\STATE $u_H$ $\Leftarrow$ DSOLVE($A_H u_H=f_H$) $\quad\quad \quad\;\;\,\rhd$ $u_H$ is FD solution of (\ref{FD}) with mesh sizes $H_x, H_y, H_z$.
\STATE $u_{H/2}$ $\Leftarrow$ DSOLVE($A_{H/2} u_{H/2}=f_{H/2}$) $\quad \rhd$ $u_{H/2}$ is FD solution of (\ref{FD}) with mesh sizes $H_x/2, H_y/2, H_z/2$.
\STATE $h_x=H_x/2, h_y=H_y/2, h_z=H_z/2$
  \FOR {$i=1$ to $L$}
   \STATE $h_x=h_x/2, h_y=h_y/2, h_z=h_z/2$
   \STATE ${w}_{h} = \textrm{EXP}_{finite}(u_{2h}, u_{4h})$  $\quad\quad\quad\quad\; \rhd$ $w_h$ is a fifth-order approximation of the actual numerical solution $u_h$, and it serves as the initial guess for the CG solver on the next finer grid.
    \WHILE {$||A_h u_h -f_h||_2>\epsilon \cdot ||f_h||_2 $}
        \STATE $u_h \Leftarrow$ CG$(A_h, u_h, f_h)$
    \ENDWHILE
    \STATE $\tilde{u}_{h} = \textrm{EXP}_{true}(u_{h}, u_{2h})$ $\quad\quad\;\;\quad\quad  \rhd$ Optional step. $\tilde{u}_{h}$ is a sixth-order approximation solution for sufficiently smooth $u$.
  \ENDFOR
\end{algorithmic}
\end{algorithm}

%
%
%
%
%
%

In Algorithm \ref{alg:oEXCMG},  the  coarsest grid  has the mesh sizes  $H_x, H_y, H_z$, the positive integer $L$ is the total number of grids except first two embedded grids, which indicates that the mesh sizes of the finest grid are $\frac{H_x}{2^{L+1}}, \frac{H_y}{2^{L+1}}, \frac{H_z}{2^{L+1}}$. DSOLVE is a direct solver used on the first two coarse grids (see line 1-2 in Algorithm \ref{alg:oEXCMG}). Procedure $\textrm{EXP}_{finite}(u_{2h}, u_{4h})$ denotes a fifth-order approximation to the actual compact FD solution $u_h$ obtained by Richardson extrapolation and tri-quartic Lagrange interpolation from the numerical solutions $u_{2h}$ and $u_{4h}$. And there is an optional step in the above algorithm  (see line 10 in Algorithm~\ref{alg:oEXCMG}), where $\textrm{EXP}_{true}(u_{h}, u_{2h})$ denotes a higher-order approximation solution on entire fine grid with mesh size $h$ from two fourth-order FD solutions $u_h$ and $u_{2h}$.  This optional step is used to increase the order of solution accuracy from fourth order to sixth order (see Table~\ref{table1}-\ref{table8} in section~\ref{sec5} for details) when the exact solution $u$ of elliptic equation (\ref{bvp}) is sufficiently smooth.

The detailed procedures of extrapolation and tri-quartic Lagrange interpolation are described in the next two subsections \ref{extra} and \ref{3D}. The differences between our new EXCMG method and existing EXCMG method~\cite{Chen1,Chen2} are listed as follows:

\begin{enumerate}[(1)]
   \item In our new EXCMG method, a fourth-order compact difference scheme, rather than  the second-order linear FE method, is employed to discretize the 3D Poisson equation.

 \item Instead of performing a fixed number of smoothing iterations as used in the existing EXCMG method~\cite{Chen1,Chen2}, a relative residual tolerance $\epsilon$ is introduced for the smoother in our EXCMG method (see line 7 in Algorithm \ref{alg:oEXCMG}), which enables us to conveniently obtain the numerical solution with the desired accuracy.

  \item    In the existing EXCMG literature~\cite{Chen1,Chen2}, a third-order approximation to the second-order FE solution is constructed to serve as the initial guess for the iterative solver on the next finer grid, and the construction of the third-order approximation to the second-order FE solution is done at every single coarse hexahedral element. However, in our new EXCMG method,  a fifth-order approximation to the fourth-order FD solution, obtained through the Richardson extrapolation and tri-quartic Lagrange interpolation, is used as the initial guess for the iterative solver. In addition, the tri-quartic interpolation should be done for every cell which  contains eight neighboring coarse hexahedral elements as shown in Fig.~\ref{Fig3D}, rather than every single coarse hexahedral element.
\end{enumerate}

\subsection{Extrapolation and quartic
interpolation: 1D case}\label{extra}

The extrapolation method is an efficient procedure for increasing the solution accuracy of many
problems in numerical analysis. Marchuk and Shaidurov \cite{Marchuk1983}
systematically studied its application in the FD method in 1983. Since then, this technique has been well demonstrated
in the framework of the FD and FE methods \cite{Neittaanmaki1987,Fobmeier1989,Han1993,Wang2010,Sun2004,Rahul2006,Munyakazi2008,Tam2000,Ma2010,Marchi2013}.

In this and next subsections, we assume that the exact solution $u$ is sufficiently smooth, and we will formally explain how to use extrapolation and  quartic interpolation techniques to obtain the fifth-order approximation $w_h$ of the fourth-order FD solution on the next finer grid, which can be regarded as another important application of the extrapolation method.  
In addition, we will also show how to construct the enhanced sixth-order accurate numerical solution $\tilde{u}_h$ for the problem (\ref{bvp}).

\subsubsection{Extrapolation for the true solution}
For simplicity, we first consider the three-levels of embedded grids $Z_i(i=0,1,2)$ with mesh sizes $h_i=h_0/2^i$ in one dimension.
Suppose $u\in H^{6}(\Omega)$, from theorem 4.1 in~\cite{Berikelashuili} (taking $m=2,s=6$) and by using the result that $H^{2}(\Omega)$ can be continuously embedded into $L^{\infty}(\Omega)$, we can get that the error $||e^i||_{\infty}$ should be $O(h^4)$, where $e^i=u^i-u$ is the error of the fourth-order compact FD solution $u^i$ with mesh size $h_i$.  Now we further assume that
the truncation error at node $x_k$ has the form
\begin{equation}\label{asy}
  e^i(x_k)=A(x_k)h_i^4+O(h_i^6),
\end{equation}
where $A(x)$ is a suitably smooth function independent of $h_i$. The truncation error expansion (\ref{asy}) will be verified by
numerical results in section \ref{sec5}.

It is well known that the traditional extrapolation is possible only at coarse grid points, where at least two approximations,
corresponding to different mesh size, are known.
From eq. (\ref{asy}), we easily obtain the Richardson extrapolation formula at coarse grid points
\begin{equation}\label{t1}
  \tilde{u}_k^1 := \frac{16 u^1_k - u^0_k}{15} = u(x_k) + O(h_0^6),\ \ k=j,j+1,
\end{equation}
which is a sixth-order approximation to the true solution at the coarse grid points.

In fact, by using the linear interpolation formula, one can also obtain a sixth-order accurate approximation at the fine grid point $x_{j+1/2}$.
Setting $i=0$ and $i=1$ in eq. (\ref{asy}) and then subtracting each other, we have
\begin{equation}\label{aaa}
  A(x_k) = \frac{16}{15h_0^4}(u_k^0-u_k^1) + O(h_0^2),\ \ k=j,j+1.
\end{equation}
From the error estimate of the linear interpolation
\begin{equation}\label{A_12}
  A(x_{j+1/2})=\frac{1}{2}(A(x_j)+A(x_{j+1}))+O(h_0^2),
\end{equation}
and substituting eq. (\ref{aaa}) into eq. (\ref{A_12}) , we get
\begin{equation}\label{Axmid}
  A(x_{j+1/2}) = \frac{8}{15h_0^4}(u_j^0-u_j^1) + \frac{8}{15h_0^4}(u_{j+1}^0-u_{j+1}^1) + O(h_0^2).
\end{equation}
Since
\begin{equation}
  u_{j+1/2}^1 = u(x_{j+1/2}) + \frac{1}{16} A(x_{j+1/2})h_0^4 + O(h_0^6),
\end{equation}
by using (\ref{Axmid}), we obtain the following midpoint extrapolation formula:
\begin{equation}\label{chenlin}
  \tilde{u}_{j+1/2}^1 := u_{j+1/2}^1 + \frac{1}{30}(u^1_j - u^0_j + u^1_{j+1} - u^0_{j+1}) = u(x_{j+1/2}) + O(h_0^6),
\end{equation}
which is a sixth-order approximation to the true solution at the fine grid point $x_{j+1/2}$.

\subsubsection{Extrapolation for the FD solution}\label{extra_solution}
In this part, we will explain, given the fourth-order FD solutions $u^0$ and $u^1$, how to use the extrapolation and high-order interpolation techniques to construct a fifth-order (to be illustrated in subsection \ref{error}) approximation $w^2$ to the FD solution $u^2$.

\begin{figure}
   \centering
   \scalebox{0.8}{\includegraphics{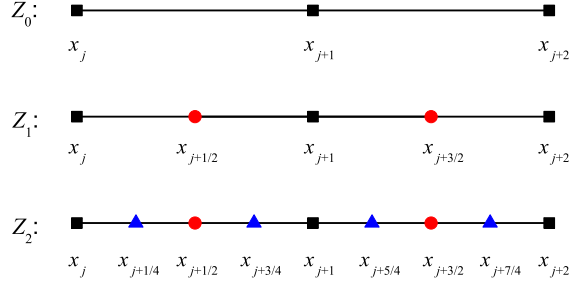}}
   \caption{Three embedded grids for two neighboring coarse elements in 1D.}\label{Fig:1}
\end{figure}

Adding one midpoint and two four equal division points,
the coarse mesh element $(x_j,x_{j+1})$ is uniformly refined into four elements of fine mesh $Z_2$ as shown in Fig. \ref{Fig:1}.
Assume there exists a constant $c$ such that
\begin{equation}\label{comb}
  c u^1+(1-c) u^0=u^2+O(h_0^6).
\end{equation} Here, we aim to use a linear combination of $u^0$ and $u^1$
to approximate the FD solution $u^2$ up to sixth-order accuracy.
Substituting the asymptotic error expansion (\ref{asy}) into (\ref{comb}), we obtain $c=17/16$ and an
extrapolation formula
\begin{equation}\label{jiedian}
   w^2_k := \frac{17 u^1_k - u^0_k}{16} = u^2_k+O(h_0^6),\ \ k=j,j+1,
\end{equation}
 at nodes $x_j$ and $x_{j+1}$. To derive the extrapolation formula at midpoint $x_{j+1/2}$,  eq. (\ref{asy}) leads to
\begin{equation}\label{bbb}
  u_{j+1/2}^2 = u_{j+1/2}^1 - \frac{15}{256} A(x_{j+1/2}) h_0^4 + O(h_0^6).
\end{equation}
Substituting eq. (\ref{Axmid}) into eq. (\ref{bbb}),  we have the following sixth-order extrapolation formula at the midpoint $x_{j+1/2}$,
\begin{align}
      w^2_{j+1/2}:=u^1_{j+1/2}+ \frac{1}{32}(u^1_j - u^0_j + u^1_{j+1} - u^0_{j+1})={u}^2_{j+1/2} + O(h_0^6). \label{zhongdian}
\end{align}
Sixth-order extrapolation formulas (\ref{jiedian}) and (\ref{zhongdian}) can be efficiently applied to each coarse-grid element $(x_j, x_{j+1})$.

Once the five approximated values $w_j^2,w_{j+1/2}^2, w_{j+1}^2,w_{j+3/2}^2$ and $w_{j+2}^2$  are obtained on the two neighboring coarse elements, we can get the following four equal division point extrapolation formulas by using the quartic interpolation
\begin{align}
      w^2_{j+1/4}&:=\displaystyle\frac1{128}\big(35w^2_{j}+140w^2_{j+1/2}-70w^2_{j+1}+28w^2_{j+3/2}-5w^2_{j+2}\big),\label{sifen1}\\
      w^2_{j+3/4}&:=\displaystyle\frac1{128}\big(-5w^2_{j}+60w^2_{j+1/2}+90w^2_{j+1}-20w^2_{j+3/2}+3w^2_{j+2}\big),\label{sifen2}\\
      w^2_{j+5/4}&:=\displaystyle\frac1{128}\big(-5w^2_{j}+28w^2_{j+1/2}-70w^2_{j+1}+140w^2_{j+3/2}+35w^2_{j+2}\big),\label{sifen3}\\
      w^2_{j+7/4}&:=\displaystyle\frac1{128}\big(3w^2_{j}-20w^2_{j+1/2}+90w^2_{j+1}+60w^2_{j+3/2}-5w^2_{j+2}\big).\label{sifen4}
\end{align}
Until now,  we have obtained a high-order approximation $w^2$ to the FD solution $u^2$, which can be used as the initial guess of the iterative solution on the fine mesh $Z_2$.

 \subsection{Extrapolation and quartic
interpolation: 3D case}\label{3D}

\begin{figure}
  \centering
  \includegraphics[width=0.9\textwidth]{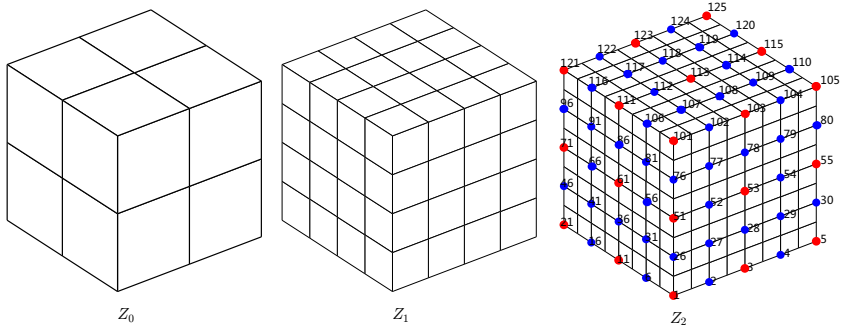}
  \caption{Three embedded grids on one interpolation cell which contains eight neighboring coarse hexahedral elements.
  }\label{Fig3D}
  \end{figure}
In this subsection, we explain how to obtain a fifth-order accurate
 approximation ${w}^2$ to the fourth-order FD solution $u^2$, and a sixth-order accurate approximate solution $\tilde{u}^1$ to the problem (\ref{bvp})
 for embedded hexahedral grids as shown in Fig. \ref{Fig3D}.

Taking every interpolation cell which consists of eight neighboring coarse hexahedral elements (see Fig. \ref{Fig3D}) into account, the construction processes of the approximation $w^2$ are as follows:
\begin{description}

  \item[Corner Nodes (such as 1, 3, 51, 53):]
   The approximate values at 27 corner nodes `\textcolor[rgb]{1,0,0}{$\medbullet$}' on such  interpolation  cell can be obtained by using the
extrapolation formula (\ref{jiedian}).

  \item[Midpoints of edges (such as 2, 6, 26, 28):]
    The approximate values at these 54 midpoints `\textcolor[rgb]{0,0,1}{$\medbullet$}' on such  interpolation  cell can be obtained by using the midpoint extrapolation formula (\ref{zhongdian}) in $x$-direction, $y$-direction or $z$-direction.

    \item[Centers of faces (such as 27, 31, 107, 109):]
    Since the center of each face on such interpolation  cell  can be viewed as the midpoint of two face diagonals, using the midpoint extrapolation formula (\ref{zhongdian})
    we can obtain two  approximate values, and take the arithmetic mean as the approximation at these 36 midpoints `\textcolor[rgb]{0,0,1}{$\medbullet$}' .

    \item[Centers of coarse hexahedral elements (such as 32, 42, 82, 92):]
        Since the center of each coarse hexahedral element on such interpolation  cell  can be viewed as the midpoint of four space diagonals, again using the midpoint extrapolation formula (\ref{zhongdian}) we can obtain four approximate values, and take the arithmetic mean as the approximation at these 8 midpoints `\textcolor[rgb]{0,0,1}{$\medbullet$}' .

\item[Other fine grid points:] The approximate values of remaining 604$(9^3-5^3)$ fine grid points on such the interpolation  cell can be obtained by using tri-quartic Lagrange interpolation with the known 125-node (27 corner nodes, 54 midpoints of edges, 36 centers of faces and 8 centers of coarse hexahedral elements) values.
\end{description}

The tri-quartic Lagrange interpolation function in terms of natural coordinates ($\xi,\eta,\zeta$) is
    \begin{equation}\label{inter}
      w^2(\xi,\eta,\zeta)=\sum_{m=1}^{125} N_m(\xi,\eta,\zeta) w_m^2,
    \end{equation}
    where the shape functions $N_m$ can be written as follows
  \begin{equation}
    {N_m}(\xi ,\eta ,\zeta ) = l_i^4(\xi)l_j^4(\eta)l_k^4(\zeta),
  \end{equation}
  where $l_i^4(x)\ (0\leq i \leq 4)$ is the Lagrange fundamental polynomials of degree 4, defined as
  \begin{equation}
    l_i^4(\xi) = \prod_{k=0, k\neq i}^{4} \frac{\xi-\xi_k}{\xi_i-\xi_k},
  \end{equation}
and $(\xi_i,\eta_j,\zeta_k)$ is the natural coordinate of node $m\; (1\leq m \leq 125)$.

When constructing the sixth-order accurate solution  $\tilde{u}^1$ based on two fourth-order accurate solutions $u^0$ and $u^1$,
the Richardson extrapolation formula (\ref{t1}) can be directly used for coarse grid points,
while the sixth-order midpoint extrapolation formula (\ref{chenlin}) can be directly used for all other fine grid points, which is
similar to the process (excluding the tri-quartic interpolation) of constructing the approximation $w^2$ described as above.

\begin{remark}
Since the compact FD solution $u_h$ of (\ref{compact}) is a fourth-order approximation of the exact solution $u$, in order to get a quite good initial guess $w_h$ for the CG solver, a tri-quartic Lagrange interpolation method is employed in this paper so that a fifth-order approximation of $w_h$ to $u_h$ is achieved. Moreover, the relative effect of how $w_h$ approximates $u_h$ becomes better when mesh is refined, thus, the number of iterations will be reduced most significantly on the finest grid, which is particularly important for solving large linear systems and can greatly reduce the computational cost. We note that the tri-quadratic interpolation used in~\cite{Pan2015} produces a third-order approximation to the second-order FE solution, and the tri-quadratic interpolation is accurate enough in that case. However, when $u_h$ is obtained from the fourth-order compact FD method  as shown in this paper, the tri-quadratic interpolation can not provide a sufficiently accurate initial guess $w_h$, the relative effect of how $w_h$ approximates $u_h$ will become worse when mesh is refined.
\end{remark}

\begin{remark}
Tri-quartic Lagrange interpolation defined by eq. (\ref{inter}) is a local operation defining on each interpolation cell containing eight neighbouring coarse elements. In fact, eq. (\ref{inter}) defines a same $(604\times 125)$  interpolation matrix on every interpolation cell, thus the approximate values of remaining 604$(9^3-5^3)$ fine-grid points on every interpolation cell  can be obtained by multiplying the  $(604\times 125)$  interpolation matrix with the vector consisting of 125 known values on such interpolation cell. Therefore, the fifth-order approximation of FD solution $w_h$ on the entire domain can be obtained  very effectively by applying the extrapolation formulas  (\ref{jiedian}) and (\ref{zhongdian}) to the 125 nodes mentioned above, and  running the  tri-quartic Lagrange interpolation (\ref{inter}) based on such 125 known values  for every interpolation cell in the entire domain.

\end{remark}


\subsection{The error analysis of initial guess $w^2$}\label{error}

Let $e=w^2-{u}^2$ be  the difference between the initial guess $w^2$ and the FD solution $u^2$. Assume that $e$ has continuous derivatives up to order 5 on interval $[x_j,x_{j+2}]$.
From (\ref{jiedian}) and (\ref{zhongdian}) we obtain the equation
 \begin{equation}\label{err}
   e(x_k) = O(h_0^6),\quad k=j,j+1/2,j+1,j+3/2,j+2.
 \end{equation}
From polynomial interpolation theory, the
error of quartic interpolation $I_4f$ can be represented as
\begin{equation}
  R_4(x) = e - I_4e = \frac{1}{5!}e^{(5)}(\xi) (x-x_j)(x-x_{j+1/2})(x-x_{j+1})(x-x_{j+3/2})(x-x_{j+2}),
\end{equation}
where $\xi \in (x_j,x_{j+2})$ depends on $x$.
Especially at four equal division points we have
\begin{equation}\label{si1}
  R_4(x_{j+1/4}) = \frac{7 h_0^5}{8\times 4^5}e^{(5)}(\xi_1) = \frac{7 h_0^5}{8192}e^{(5)}(x_{j+1}) + o(h_0^5),
\end{equation}
\begin{equation}\label{si2}
  R_4(x_{j+3/4}) = -\frac{3 h_0^5}{8\times 4^5}e^{(5)}(\xi_2) = -\frac{3 h_0^5}{8192}e^{(5)}(x_{j+1}) + o(h_0^5),
\end{equation}
and
\begin{equation}\label{si3}
  R_4(x_{j+5/4}) = \frac{3 h_0^5}{8\times 4^5}e^{(5)}(\xi_3) = \frac{3 h_0^5}{8192}e^{(5)}(x_{j+1}) + o(h_0^5) \approx -R_4(x_{j+3/4}),
\end{equation}
\begin{equation}\label{si4}
  R_4(x_{j+7/4}) = -\frac{7 h_0^5}{8\times 4^5}e^{(5)}(\xi_4) = -\frac{7 h_0^5}{8192}e^{(5)}(x_{j+1}) + o(h_0^5)\approx -R_4(x_{j+1/4}).
\end{equation}

It follows from eqs. (\ref{err}) and  (\ref{si1})-(\ref{si4}) that
\begin{equation}\label{err_si}
   e(x_k) = I_4 e(x_k) + R_4(x_k) = O(h_0^5) ,\quad k=j+1/4,j+3/4,j+5/4,j+7/4,
 \end{equation}
 which means that the initial guess $w^2$ obtained by extrapolation and  quartic interpolation is a fifth-order accurate approximation to the FD solution $u^2$.

The above error analysis can be directly extended to 3D case (see numerical verification in Section \ref{sec5}: the last columns in Table \ref{exam1}-\ref{table9}).
In addition, eqs.(\ref{si3}) and (\ref{si4}) imply that the initial error $e(x)$ forms a high-frequency oscillation in the entire domain,
however, it can be smoothed out after a few CG iterations (see Fig. \ref{RRe} for details).

\section{Numerical experiments}\label{sec5}
In this section, in order to illustrate the  efficiency  of the new EXCMG method comparing to the classical V-cycle and W-cycle multigrid methods with  the
Gauss-Seidel relaxation and the CG relaxation, we present the numerical results for six examples with smooth and finite regular solutions using the proposed method. Our code is written in Fortran 90 with double precision arithmetic, and compiled with  Intel Visual Fortran Compiler XE 12.1 under 64-bit Windows 7. All programs are carried out on a personal desktop equipped with Intel(R) Core(TM) i7-4790K CPU (4.00 GHz) and 16GB RAM.

The order of convergence of the method is computed by
\begin{equation}
  \textrm{order} = \log_2\frac{||u_h-u||}{||u_{h/2}-u||},
\end{equation}
where $||\cdot||$ denotes some norm (for instance, $L^2$-norm or $L^\infty$-norm) and $u$ is the true solution.

\subsection{Numerical accuracy}

\begin{example}\label{exam1}
  The test Problem \ref{exam1} can be written as
\begin{equation}\label{prob1}
         \frac{\partial^2 u}{\partial x^2}+\frac{\partial^2 u}{\partial y^2}+\frac{\partial^2 u}{\partial z^2} =e^{z}\sin(xy)(1-x^2-y^2),\quad  \textrm{in } \Omega=[0,1]^3,
\end{equation}
where the boundary conditions are
\begin{align*}
& u(0,y,z)=u(x,0,z)=0,\quad u(1,y,z)=e^z\sin(y),\quad u(x,1,z)=e^z\sin(x),
\end{align*}
and
\begin{align*}
u(x,y,0)=\sin(xy),\quad u(x,y,1)=e\sin(xy).
\end{align*}
The analytic solution of eq. (\ref{prob1}) is
$$u(x,y,z)=e^z\sin(xy),$$
which is a sufficiently smooth function.
\end{example}

Using 7 embedded grids with the coarsest grid $4\times4\times4$, we present the numerical results for Problem \ref{exam1} obtained by the new EXCMG method with $\epsilon=10^{-14}$ in Table~\ref{table1}-\ref{table2}.  Table~\ref{table1} lists  the $L^2$-error of the compact FD solution $u_h$, the $L^2$-error of the gradient of the FD solution $\nabla u_h$, the $L^2$-error {of} the extrapolated solution $\tilde{u}_h$, the $L^2$-norm of the difference between the initial guess $w_h$ and the FD solution $u_h$, and corresponding  convergence rates.
 Table~\ref{table2} gives all errors and convergence rates in $L^{\infty}$-norm. Since a direct solver is used for the first two coarse levels of grids, we only list the results starting from the third level of grid $16\times16\times16$.

Here we explain how to numerically compute the gradient $\nabla u_h$ after we obtain the FD solution $u_h$. First, we use the following
fourth-order, one-sided, FD approximation of the partial derivative $u_x$ on the boundary grid points,
\begin{align*}
&(u_x)_{0,j,k}=-\frac{25}{12h_x}u_{0,j,k}+\frac{4}{h_x}u_{1,j,k}-\frac{3}{h_x}u_{2,j,k}+\frac{4}{3h_x}u_{3,j,k}-\frac{1}{4h_x}u_{4,j,k},\ \textrm{for}\ j=0,\cdots,N_y, k=0,\cdots, N_z,\\
&(u_x)_{N_x,j,k}=\frac{25}{12h_x}u_{N_x,j,k}-\frac{4}{h_x}u_{N_x-1,j,k}+\frac{3}{h_x}u_{N_x-2,j,k}-\frac{4}{3h_x}u_{N_x-3,j,k}+\frac{1}{4h_x}u_{N_x-4,j,k}, \ \textrm{for}\ j=0,\cdots,N_y, k=0,\cdots, N_z.
\end{align*}
Then we can obtain $(u_x)_{i,j,k}, (i=1,\cdots, N_x-1)$ on the internal grid points by solving the following linear system resulting from the fourth-order compact FD scheme~\cite{Collatz},
\begin{align*}
\frac{1}{6}(u_x)_{i-1,j,k}+\frac{4}{6}(u_x)_{i,j,k}+\frac{1}{6}(u_x)_{i+1,j,k}=\frac{u_{i+1,j,k}-u_{i-1,j,k}}{2h_x}, \ \ \textrm{for}\ j=0,\cdots,N_y, k=0,\cdots, N_z.
\end{align*}
The above 1D tridiagonal system can be solved fast by the Thomas algorithm. Clearly, we can get $u_y$ and $u_z$ from similar procedures.
Then, $\nabla u_h$ can be obtained efficiently.

As we can see from table \ref{table1}-\ref{table2} that initial guess $w_h$ is a fifth-order approximation to the FD solution $u_h$, which validates our theoretical analysis in section \ref{error}, and the FD solution $u_h$ achieves the full fourth-order accuracy. The numerical gradient $\nabla u_h$ is also a fourth-order approximation to the exact gradient $\nabla u$ in both the $L^2$-norm and $L^{\infty}$-norm, while the extrapolated solution $\tilde{u}_h$ converges with sixth-order accuracy on all grids except the finest grid.
This is due to the fact that the extrapolated solution $\tilde{u}_h$ is obtained from two fourth-order FD solutions $u_h$ and $u_{2h}$,
these two solutions must be extremely accurate in order to obtain a sixth-order accurate solution $\tilde{u}_h$.
As the grid becomes finer, the relative residual tolerance needs to be smaller. Thus, the extrapolated solution $\tilde{u}_h$ starts to lose convergence order when the grid is fine enough since a uniform tolerance is used in our EXCMG algorithm.
And in this example, on the finest mesh $256\times256\times256$, the maximum error between the extrapolated solution $\tilde{u}_h$ and the exact solution $u$ already reaches $O(10^{-14})$, which is almost the machine accuracy, although the method does not achieve the full sixth-order on the finest grid. Additionally, we can see that the numerical results confirm with the asymptotic  error  expansion (\ref{asy}).

\begin{table}[!tbp]
\tabcolsep=7pt
\caption{Errors and convergence rates with $\epsilon=10^{-14}$ in $L^{2}$-norm for Example \ref{exam1}.} \centering
\begin{tabular}{llrlrlrlr}
    \hline
     {mesh}& $||u_h-u||_{2}$ &   order & $||\nabla(u_h-u)||_{2}$ &   order &   $||\tilde{u}_h-u||_{2}$ &   order  &   $||w_h-u_h||_{2}$ &   order  \\
\hline
$  16\times  16\times  16$ &  $1.67(-08)$ &      & $1.08(-06)$ &      & $1.38(-09) $  &       &  $4.36(-07)$ &      \\
$  32\times  32\times  32$ &  $1.09(-09)$ & 3.93 & $4.74(-08)$ & 4.51 & $2.40(-11) $  & 5.84  &  $1.29(-08)$ & 5.08  \\
$ 64\times  64\times  64$  &  $7.00(-11)$ & 3.97 & $2.12(-09)$ & 4.48 & $3.91(-13) $  & 5.94  &  $3.94(-10)$ & 5.04  \\
$ 128\times 128\times 128$ &  $4.42(-12)$ & 3.98 & $9.70(-11)$ & 4.45 & $6.22(-15) $  & 5.98  &  $1.22(-11)$ & 5.02    \\
$ 256\times 256\times 256$ &  $2.82(-13)$ & 3.97 & $4.62(-12)$ & 4.39 & $5.68(-15) $  & 0.13  &  $3.79(-13)$ & 5.01    \\
\hline
\end{tabular}\label{table1}
\end{table}

\begin{table}[!tbp]
\tabcolsep=7pt
\caption{Errors and convergence rates with $\epsilon=10^{-14}$ in $L^{\infty}$-norm for Example \ref{exam1}.} \centering
\begin{tabular}{llrlrlrlr}
    \hline
     {mesh}& $||u_h-u||_{\infty}$ &   order & $||\nabla(u_h-u)||_{\infty}$ &   order &  $||\tilde{u}_h-u||_{\infty}$ &   order  &   $||w_h-u_h||_{\infty}$ &   order  \\
\hline
$  16\times  16\times  16$ & $5.47(-08)$  &      &  $8.25(-06)$  &      & $9.12(-09) $  &       &  $3.85(-06)$  &        \\
$  32\times  32\times  32$ & $3.43(-09)$  & 4.00 &  $5.18(-07)$  & 3.99 & $1.77(-10) $ & 5.68  &  $1.31(-07)$  & 4.88   \\
$ 64\times  64\times   64$ & $2.15(-10)$  & 4.00 &  $3.24(-08)$  & 4.00 & $3.20(-12)$  & 5.79  &  $4.19(-09)$  & 4.96   \\
$ 128\times 128\times 128$ & $1.34(-11)$  & 4.00 &  $2.03(-09)$  & 4.00 & $5.42(-14)$  & 5.88  &  $1.32(-10)$ & 4.99  \\
$ 256\times 256\times 256$ & $8.49(-13)$  & 4.00 &  $1.27(-10)$  & 4.00 & $4.00(-14)$  & 0.44  &  $4.15(-12)$ & 4.99  \\
\hline
\end{tabular}\label{table2}
\end{table}

\begin{example}\label{exam2}
  The test Problem \ref{exam2} can be written as
\begin{equation}\label{prob2}
         \frac{\partial^2 u}{\partial x^2}+\frac{\partial^2 u}{\partial y^2}+\frac{\partial^2 u}{\partial z^2} = 0, \quad  \textrm{in } \Omega=[0,1]^3,
\end{equation}
where the boundary conditions are
\begin{equation*}
u(0,y,z)=e^y\sin(\sqrt{2}\,z),\quad u(x,0,z)=e^x\sin(\sqrt{2}\,z),\quad u(x,y,0)=0,
\end{equation*}
and
\begin{equation*}
u(1,y,z)=e^{1+y}\sin(\sqrt{2}\,z),\quad u(x,1,z)=e^{x+1}\sin(\sqrt{2}\,z),\quad u(x,y,1)=e^{x+y}\sin(\sqrt{2}).
\end{equation*}
The analytic solution of eq. (\ref{prob2}) is
$$u=e^{x+y}\sin(\sqrt{2}\,z),$$
which is a harmonic function and has arbitrary order smooth derivatives.
\end{example}

Again, we use 7 embedded grids with the coarsest grid $4\times4\times4$, and the corresponding numerical results obtained by the EXCMG method with $\epsilon=10^{-14}$ are listed in table~\ref{table3} and~\ref{table4}. Once again, initial  guess $w_h$ is a fifth-order approximation of the FD solution $u_h$, the FD solution $u_h$ is fourth-order accurate, and the numerical gradient $\nabla u_h$ is also a fourth-order approximation to the exact gradient $\nabla u$,
while the extrapolated solution $\tilde{u}_h$ converges to exact solution $u$ with  sixth-order but starts to lose accuracy on the finest grid $256\times 256\times 256$. Additionally, the maximum error between the extrapolated solution $\tilde{u}_h$ and the exact solution $u$ is less than  $6.0\times 10^{-14}$, which means that the solution $\tilde{u}_h$ is already accurate enough, and we don't need to further reduce the error tolerance.

\begin{table}[!tbp]
\tabcolsep=7pt
\caption{Errors and convergence rates with $\epsilon=10^{-14}$ in $L^{2}$-norm for Example \ref{exam2}.} \centering
\begin{tabular}{llrlrlrlr}
    \hline
     {mesh}& $||u_h-u||_{2}$ &   order & $||\nabla(u_h-u)||_{2}$ &   order &   $||\tilde{u}_h-u||_{2}$ &   order  &   $||w_h-u_h||_{2}$ &   order  \\
\hline
$  16\times  16\times  16$  & $4.26(-08)$ &      & $150(-05)$ &      & $2.28(-09) $  &       &  $6.21(-06)$ &     \\
$  32\times  32\times  32$  & $2.79(-09)$ & 3.94 & $663(-07)$ & 4.50 & $3.88(-11) $  & 5.88  &  $1.95(-07)$ &5.00 \\
$ 64\times  64\times  64$   & $1.78(-10)$ & 3.97 & $295(-08)$ & 4.49 & $6.28(-13) $  & 5.95  &  $6.10(-09)$ &5.00 \\
$ 128\times 128\times 128$  & $1.13(-11)$ & 3.98 & $132(-09)$ & 4.48 & $1.02(-14)$  & 5.94  &   $1.91(-10)$ &5.00    \\
$ 256\times 256\times 256$  & $7.32(-13)$ & 3.94 & $596(-11)$ & 4.47 & $3.11(-14)$  & $-$1.60  &  $5.97(-12)$ &5.00    \\
\hline
\end{tabular}\label{table3}
\end{table}

\begin{table}[!tbp]
\tabcolsep=7pt
\caption{Errors and convergence rates with $\epsilon=10^{-14}$ in $L^{\infty}$-norm for Example \ref{exam2}.} \centering
\begin{tabular}{llrlrlrlr}
    \hline
     {mesh} &$||u_h-u||_{\infty}$ &   order & $||\nabla(u_h-u)||_{\infty}$ &   order &   $||\tilde{u}_h-u||_{\infty}$ &   order  &   $||w_h-u_h||_{\infty}$ &   order  \\
\hline
$  16\times  16\times  16$  &$1.16(-07)$ &      &$1.26(-4)$ &     &  $1.01(-08) $ &       &  $3.39(-05)$ &        \\
$  32\times  32\times  32$  &$7.22(-09)$ & 4.00 &$7.95(-6)$ & 3.99&  $1.86(-10) $ &  5.76 &  $1.11(-06)$ & 4.94   \\
$ 64\times 64\times 64$     &$4.52(-10)$ & 4.00 &$4.98(-7)$ & 4.00&  $3.17(-12)$  &  5.87 &  $3.54(-08)$ & 4.96   \\
$ 128\times 128\times 128$  &$2.83(-11)$ & 4.00 &$3.11(-8)$ & 4.00&  $5.95(-14)$  & 5.74  & $1.12(-09)$ & 4.99  \\
$ 256\times 256\times 256$  &$1.81(-12)$ & 3.96 &$1.95(-9)$ & 4.00&  $1.07(-13)$  & $-$0.85  & $3.51(-11)$ & 4.99  \\
\hline
\end{tabular}\label{table4}
\end{table}

\begin{example}\label{gaussian}
  The test Problem \ref{gaussian} can be written as
\begin{equation}\label{prob5}
\left\{ \begin{aligned}
         \frac{\partial^2 u}{\partial x^2}+\frac{\partial^2 u}{\partial y^2}+\frac{\partial^2 u}{\partial z^2} &=f(x,y,z),  &\textrm{in } &\Omega=[0,1]^3,\\
          u&=g(x,y,z),  &\textrm{on } &\partial \Omega,
        \end{aligned} \right.
\end{equation}
where $f$ and $g$ are determined from the exact solution
$$u=e^{-3\big((x-0.5)^2+(y-0.5)^2+(z-0.5)^2\big)},$$
which is a 3D Gaussian function. It varies rapidly near the point $(0.5,0.5,0.5)$.
\end{example}

Table~\ref{table31} and~\ref{table32} list the numerical results obtained by the EXCMG method with $\epsilon=10^{-11}$. One more time, one can see that initial  guess $w_h$ is a fifth-order approximation of the FD solution $u_h$, the FD solution $u_h$ is fourth-order accurate (although the convergent order is slightly reduced on the finest grid), and the numerical gradient $\nabla u_h$ is also a fourth-order approximation to the exact gradient $\nabla u$,
while the extrapolated solution $\tilde{u}_h$  is  sixth-order accurate. Therefore, our EXCMG method is still very effective for the problem with very rapid variations.

\begin{table}[!tbp]
\tabcolsep=7pt
\caption{Errors and convergence rates with $\epsilon=10^{-11}$ in $L^{2}$-norm for Example \ref{gaussian}.} \centering
\begin{tabular}{llrlrlrlr}
    \hline
     {mesh}& $||u_h-u||_{2}$ &   order & $||\nabla(u_h-u)||_{2}$ &   order &   $||\tilde{u}_h-u||_{2}$ &   order  &   $||w_h-u_h||_{2}$ &   order  \\
\hline
$  16\times  16\times  16$ &  $1.22(-06)$ &      & $4.08(-04)$ &      & $6.68(-08) $  &       &  $2.29(-04)$ &      \\
$  32\times  32\times  32$ &  $7.90(-08)$ & 3.95 & $1.29(-05)$ & 4.99 & $1.15(-09) $  & 5.86  &  $5.80(-06)$ & 5.30  \\
$ 64\times  64\times  64$  &  $5.04(-09)$ & 3.97 & $4.42(-07)$ & 4.86 & $1.87(-11) $  & 5.94  &  $1.86(-07)$ & 4.96  \\
$ 128\times 128\times 128$ &  $3.19(-10)$ & 3.98 & $1.73(-08)$ & 4.67 & $3.05(-13) $  & 5.94  &  $5.86(-09)$ & 4.99    \\
$ 256\times 256\times 256$ &  $2.40(-11)$ & 3.73 & $7.50(-10)$ & 4.53 & $4.52(-12) $  &$-$3.89&  $1.84(-10)$ & 4.99    \\
\hline
\end{tabular}\label{table31}
\end{table}

\begin{table}[!tbp]
\tabcolsep=7pt
\caption{Errors and convergence rates with $\epsilon=10^{-11}$ in $L^{\infty}$-norm for Example \ref{gaussian}.} \centering
\begin{tabular}{llrlrlrlr}
    \hline
     {mesh}& $||u_h-u||_{\infty}$ &   order & $||\nabla(u_h-u)||_{\infty}$ &   order &  $||\tilde{u}_h-u||_{\infty}$ &   order  &   $||w_h-u_h||_{\infty}$ &   order  \\
\hline
$  16\times  16\times  16$ & $4.80(-06)$  &      &  $1.03(-03)$  &      & $2.14(-07) $  &       &  $1.13(-03)$  &        \\
$  32\times  32\times  32$ & $2.97(-07)$  & 4.01 &  $4.35(-05)$  & 4.56 & $4.07(-09) $  & 5.71  &  $2.50(-05)$  & 5.50   \\
$ 64\times  64\times   64$ & $1.85(-08)$  & 4.00 &  $2.01(-06)$  & 4.43 & $6.63(-11)$   & 5.94  &  $9.98(-07)$  & 4.65   \\
$ 128\times 128\times 128$ & $1.16(-09)$  & 4.00 &  $1.03(-07)$  & 4.28 & $1.04(-12)$   & 5.99  &  $3.02(-08)$  & 5.05  \\
$ 256\times 256\times 256$ & $8.95(-11)$  & 3.69 &  $5.55(-09)$  & 4.22 & $1.85(-11)$   &$-$4.15&  $9.27(-10)$  & 5.03  \\
\hline
\end{tabular}\label{table32}
\end{table}

\begin{example}\label{exam4}
  The test Problem \ref{exam4} can be written as
\begin{equation}\label{prob4}
         \frac{\partial^2 u}{\partial x^2}+\frac{\partial^2 u}{\partial y^2}+\frac{\partial^2 u}{\partial z^2} =-5.25 \pi^2 \sin(2\pi x)\sin(\pi y)\sin(\frac{\pi}{2}z),\quad  \textrm{in } \Omega=[0,1]^3,
\end{equation}
where the boundary conditions are
\begin{equation*}
u(0,y,z)=u(1,y,z)=u(x,0,z)=u(x,1,z)=u(x,y,0)=0\;\; \textrm{and } u(x,y,1)=\sin(2\pi x)\sin(\pi y).
\end{equation*}
The analytic solution of eq. (\ref{prob4}) is
$$u(x,y,z)=\sin(2\pi x)\sin(\pi y)\sin(\frac{\pi}{2}z),$$
which changes fastest in the $x$ direction, faster in the $y$ direction and slowest in the $z$ direction.
\end{example}


Since the solution has the fastest change in the $x$-direction and the slowest change in the $z$-direction, we use the coarsest grid $8\times4\times2$
in the EXCMG algorithm. Table~\ref{table5} and~\ref{table6} list the numerical data  obtained by EXCMG method using a  tolerance $\epsilon=10^{-9}$.
Again, the initial  guess $w_h$ is a fifth-order approximation of the FD solution $u_h$, the FD solution $u_h$ is fourth-order accurate, and the numerical gradient $\nabla u_h$ is also a fourth-order approximation to the exact gradient $\nabla u$, while the extrapolated solution $\tilde{u}_h$ achieves sixth-order accuracy but starts to lose accuracy on the  finest grid since a uniform tolerance $\epsilon=10^{-9}$ is used on each level of grid. The maximum error between the extrapolated solution $\tilde{u}_h$ and the exact solution $u$ already reaches $O(10^{-11})$ on the finest grid which is again quite accurate.

\begin{table}[!tbp]
\tabcolsep=7pt
\caption{Errors and convergence rates with $\epsilon=10^{-9}$ in $L^{2}$-norm for Example \ref{exam4}.} \centering
\begin{tabular}{llrlrlrlr}
    \hline
     {mesh}& $||u_h-u||_{2}$ &   order & $||\nabla(u_h-u)||_{2}$ &   order &   $||\tilde{u}_h-u||_{2}$ &   order  &   $||w_h-u_h||_{2}$ &   order  \\
\hline
$  32\times 16\times  8$ & $3.58(-06)$  &     & $2.85(-4)$  &     & $4.55(-07) $  &        &  $7.86(-05)$ &        \\
$  64\times 32\times 16$ & $2.35(-07)$  & 3.93& $1.37(-5)$  & 4.38& $8.36(-09) $  &  5.76  &  $2.58(-06)$ & 4.93   \\
$ 128\times 64\times 32$ & $1.51(-08)$  & 3.96& $6.34(-7)$  & 4.43& $1.39(-10)$   &  5.92  &  $8.19(-08)$ & 4.98   \\
$ 256\times128\times 64$ & $9.51(-10)$  & 3.98& $2.96(-8)$  & 4.42& $3.24(-12)$   &  5.42  &  $2.58(-09)$ & 4.99  \\
$ 512\times256\times128$ & $5.73(-11)$  & 4.05& $1.67(-9)$  & 4.15& $9.12(-12)$   & $-$1.49  &  $7.99(-11)$ & 5.01  \\
\hline
\end{tabular}\label{table5}
\end{table}

\begin{table}[!tbp]
\tabcolsep=7pt
\caption{Errors and convergence rates with $\epsilon=10^{-9}$ in $L^{\infty}$-norm for Example \ref{exam4}.} \centering
\begin{tabular}{llrlrlrlr}
    \hline
     {mesh}& $||u_h-u||_{\infty}$ &   order & $||\nabla(u_h-u)||_{\infty}$ &   order &  $||\tilde{u}_h-u||_{\infty}$ &   order  &   $||w_h-u_h||_{\infty}$ &   order  \\
\hline
$  32\times 16\times  8$ & $1.10(-05)$  &     & $1.78(-3)$  &     & $2.74(-06) $  &        &  $3.58(-04)$ &        \\
$  64\times 32\times 16$ & $6.97(-07)$  & 3.98& $1.16(-4)$  & 3.93& $6.23(-08) $  &  5.46  &  $9.97(-06)$ & 5.17   \\
$ 128\times 64\times 32$ & $4.35(-08)$  & 4.00& $7.36(-6)$  & 3.98& $1.18(-09)$   &  5.72  &  $3.20(-07)$ & 4.96   \\
$ 256\times128\times 64$ & $2.71(-09)$  & 4.01& $4.61(-7)$  & 4.00& $1.92(-11)$   &  5.94  &  $9.59(-09)$ & 5.06  \\
$ 512\times256\times128$ & $1.82(-10)$  & 3.90& $3.74(-8)$  & 3.63& $4.58(-11)$   & $-$1.26  &  $2.82(-10)$ & 5.09  \\
\hline
\end{tabular}\label{table6}
\end{table}

Previous examples are results for the 3D Poisson equation where the exact solution is infinitely many times continuously differentiable.
In the following examples, we will show the results using the new EXCMG method for the cases where the exact
solutions have finite regularities.

\begin{example}\label{exam5}
  The test Problem \ref{exam5} can be written as
\begin{equation}\label{prob5}
\left\{ \begin{aligned}
         \frac{\partial^2 u}{\partial x^2}+\frac{\partial^2 u}{\partial y^2}+\frac{\partial^2 u}{\partial z^2} &=f(x,y,z),  &\textrm{in } &\Omega=[0,1]^3,\\
          u&=g(x,y,z),  &\textrm{on } &\partial \Omega,
        \end{aligned} \right.
\end{equation}
where $f(x,y,z)$ and  $g(x,y,z)$ are determined from the exact solution
$$u=\frac{x^3y^3z^3}{(x^2+y^2+z^2)^{1.5}},$$
which has singularity at the origin and belongs to $H^{7.5-\varepsilon}$ ($\varepsilon$ is an arbitrary positive constant). It follows from the Sobolev embedding theorem that $u\in C^{6-\epsilon}$.
\end{example}

\begin{table}[!tbp]
\tabcolsep=7pt
\caption{Errors and convergence rates with $\epsilon=10^{-12}$ in $L^{2}$-norm for Example \ref{exam5}.} \centering
\begin{tabular}{llrlrlrlr}
    \hline
     {mesh}& $||u_h-u||_{2}$ &   order  & $||\nabla(u_h-u)||_{2}$ &   order &    $||\tilde{u}_h-u||_{2}$ &   order  &   $||w_h-u_h||_{2}$ &   order  \\
\hline
$  16\times  16\times  16$ & $5.44(-08)$  &      & $1.22(-05)$  &      &   $3.35(-09) $ &        &   $2.32(-06)$ &     \\
$  32\times  32\times  32$ & $3.56(-09)$  & 3.94 & $6.84(-07)$  & 4.15 &   $5.37(-11) $ &   5.96 &   $1.29(-07)$ & 4.18  \\
$ 64\times  64\times  64$  & $2.27(-10)$  & 3.97 & $3.25(-08)$  & 4.39 &   $8.57(-13)$  &   5.97  &  $4.19(-09)$ & 4.94  \\
$ 128\times 128\times 128$ & $1.44(-11)$  & 3.98 & $1.47(-09)$  & 4.46 &   $1.34(-14)$  &   6.00  &  $1.32(-10)$ & 4.99    \\
$ 256\times 256\times 256$ & $9.57(-13)$  & 3.91 & $6.59(-11)$  & 4.48 &   $1.18(-13)$  & $-3.14$ &  $4.12(-12)$ & 5.00    \\
\hline
\end{tabular}\label{table7}
\end{table}

\begin{table}[!tbp]
\tabcolsep=7pt
\caption{Errors and convergence rates with $\epsilon=10^{-12}$ in $L^{\infty}$-norm for Example \ref{exam5}.} \centering
\begin{tabular}{llrlrlrlr}
    \hline
     {mesh}& $||u_h-u||_{\infty}$ &   order  & $||\nabla(u_h-u)||_{\infty}$ &   order &   $||\tilde{u}_h-u||_{\infty}$ &   order  &   $||w_h-u_h||_{\infty}$ &   order  \\
\hline
$  16\times  16\times  16$ & $1.13(-07)$  &      & $8.26(-5)$  &      & $1.18(-08) $  &        &  $2.53(-05)$ &        \\
$  32\times  32\times  32$ & $7.16(-09)$  & 3.98 & $5.83(-6)$  & 3.82 & $2.08(-10) $  & 5.83   &  $6.96(-07)$ & 5.18   \\
$ 64 \times 64 \times 64$  & $4.48(-10)$  & 4.00 & $3.76(-7)$  & 3.96 & $3.32(-12)$   & 5.97   &  $2.48(-08)$ & 4.81   \\
$ 128\times 128\times 128$ & $2.80(-11)$  & 4.00 & $2.36(-8)$  & 3.99 & $5.25(-14)$   & 5.98   &  $8.02(-10)$ & 4.95  \\
$ 256\times 256\times 256$ & $1.81(-12)$  & 3.95 & $1.50(-9)$  & 3.97 & $2.58(-13)$   & $-$2.30&  $2.53(-11)$ & 4.99  \\
\hline
\end{tabular}\label{table8}
\end{table}

In the numerical computation, we also use 7 embedded grids with the coarsest grid $4\times4\times4$, and the corresponding numerical results by the EXCMG method with $\epsilon=10^{-12}$ are listed in table~\ref{table7} and~\ref{table8}. From  table~\ref{table7} and~\ref{table8}, one can easily find that the results are the same as previous examples, i.e.,  in both $L^2$ and $L^{\infty}$-norms, the initial  guess $w_h$ is a fifth-order approximation of the FD solution $u_h$, the FD solution $u_h$ is fourth-order accurate, and the numerical gradient $\nabla u_h$ is also a fourth-order approximation to the exact gradient $\nabla u$, while the extrapolated solution $\tilde{u}_h$ achieves sixth-order accuracy but starts to lose accuracy on the  finest grid since a uniform tolerance $\epsilon=10^{-12}$ is used on each level of grid.

We further carry out the computations for other cases when the exact solution $u$ has lower regularities, we find that if the exact solution $u\in H^{s}$ ($s<7.5$), then the extrapolated solution $\tilde{u}_h$ will not  reach sixth-order accuracy in $L^{\infty}$-norm.  In addition, we find that only when the exact solution $u$ satisfies that $u\in H^{s}$ ($s\geq5.5$),  then the numerical solution $u_h$ can reach fourth-order accuracy in  $L^{\infty}$-norm. This is not surprising since $H^{7.5+\epsilon}(\Omega)$ can be continuously embedding into $C^{6}(\Omega)$ and $H^{5.5+\epsilon}(\Omega)$ can be continuously embedding into $C^{4}(\Omega)$ from the Sobolev embedding theorem.

In the final part of this section, we will show the results for one example where the exact solution $u \in H^{5.5-\varepsilon}$ ($\varepsilon$ is an arbitrary small positive constant).
\begin{example}\label{exam3}
  The test Problem \ref{exam3} can be written as
\begin{equation}\label{prob3}
\left\{ \begin{aligned}
         \frac{\partial^2 u}{\partial x^2}+\frac{\partial^2 u}{\partial y^2}+\frac{\partial^2 u}{\partial z^2} &=\frac{8xyz}{({x^2+y^2+z^2})^{0.5}},  &\textrm{in } &\Omega=[0,1]^3,\\
          u&=g(x,y,z),  &\textrm{on } &\partial \Omega,
        \end{aligned} \right.
\end{equation}
where eq. (\ref{prob3}) has singularity at the origin and $g(x,y,z)$ is determined from the exact solution
$$u=xyz(x^2+y^2+z^2)^{0.5},$$
which belongs to $H^{5.5-\varepsilon}$ ($\varepsilon$ is an arbitrary small positive constant). It follows from the Sobolev embedding theorem that $u\in C^{4-\epsilon}$.
\end{example}

Once again, we use 7 embedded grids with the coarsest grid $4\times4\times4$, and the corresponding numerical results by the EXCMG method with $\epsilon=10^{-13}$ are listed in table~\ref{table9} and~\ref{table10}. Since in this case, the exact solution $u$ is only has a finite regularity $H^{5.5-\varepsilon}$ ($\varepsilon$ is any positive constant). From table~\ref{table9} and~\ref{table10}, we can see that the numerical solution $u_h$ is a fourth-order approximation to the exact solution in both $L^2$ and $L^{\infty}$-norms. However, due to the lack of regularity of the exact solution, we can see from table~\ref{table9} and~\ref{table10}  that the numerical gradient $\nabla u_h$ converges with fourth-order accuracy in $L^2$-norm but only third-order in $L^{\infty}$-norm, the extrapolated solution $\tilde{u}_h$ is $5.5$th-order accurate in $L^2$-norm but only $3.5$th-order accurate in $L^{\infty}$-norm,  while the initial  guess $w_h$ is only a fourth-order approximation to the FD solution $u_h$ in $L^{\infty}$-norm. Nonetheless, the initial guess $w_h$ is still a fifth-order approximation to the FD solution $u_h$ in $L^2$-norm.  Since the relative residual in the CG solver in our new EXCMG method is calculated based on the $L^2$-norm (see line 7 of the algorithm~\ref{alg:oEXCMG}), thus, our EXCMG method is still effective for such low regularity  problems ($u\in H^{5.5-\varepsilon}$), and extrapolation can also help us to increase the accuracy of initial guess $w_h$ in $L^2$-norm, which would widen the scope of applicability of our method.

\begin{table}[!tbp]
\tabcolsep=7pt
\caption{Errors and convergence rates with $\epsilon=10^{-13}$ in $L^{2}$-norm for Example \ref{exam3}.} \centering
\begin{tabular}{llrlrlrlr}
    \hline
     {mesh}& $||u_h-u||_{2}$ &   order  & $||\nabla(u_h-u)||_{2}$ &   order &    $||\tilde{u}_h-u||_{2}$ &   order  &   $||w_h-u_h||_{2}$ &   order  \\
\hline
$  16\times  16\times  16$ & $2.39(-08)$  &      & $7.11(-06)$  &      &   $2.68(-09) $ &       &  $1.54(-06)$ &     \\
$  32\times  32\times  32$ & $1.55(-09)$  & 3.94 & $4.45(-07)$  & 4.00 &   $6.35(-11) $ &  5.40 &  $6.71(-08)$ & 4.52  \\
$ 64\times  64\times  64$  & $9.90(-11)$  & 3.97 & $2.46(-08)$  & 4.18 &   $1.45(-12)$  &  5.45 &  $2.39(-09)$ & 4.81  \\
$ 128\times 128\times 128$ & $6.26(-12)$  & 3.98 & $1.28(-09)$  & 4.27 &   $3.26(-14)$  &  5.48 &  $7.90(-11)$ & 4.92    \\
$ 256\times 256\times 256$ & $3.84(-13)$  & 4.03 & $6.29(-11)$  & 4.34 &   $1.03(-13)$  &$-$1.66&  $2.54(-12)$ & 4.96    \\
\hline
\end{tabular}\label{table9}
\end{table}

\begin{table}[!tbp]
\tabcolsep=7pt
\caption{Errors and convergence rates with $\epsilon=10^{-13}$ in $L^{\infty}$-norm for Example \ref{exam3}.} \centering
\begin{tabular}{llrlrlrlr}
    \hline
     {mesh}& $||u_h-u||_{\infty}$ &   order  & $||\nabla(u_h-u)||_{\infty}$ &   order &   $||\tilde{u}_h-u||_{\infty}$ &   order  &   $||w_h-u_h||_{\infty}$ &   order  \\
\hline
$  16\times  16\times  16$ & $1.25(-07)$  &      & $3.26(-5)$  &      & $1.09(-07) $  &       &  $7.05(-06)$ &        \\
$  32\times  32\times  32$ & $7.81(-09)$  & 4.00 & $4.08(-6)$  & 3.00 & $6.82(-09) $  & 3.50  &  $4.40(-07)$ & 4.00   \\
$ 64 \times 64 \times 64$  & $4.88(-10)$  & 4.00 & $5.10(-7)$  & 3.00 & $4.26(-10)$   & 3.50  &  $2.75(-08)$ & 4.00   \\
$ 128\times 128\times 128$ & $3.05(-11)$  & 4.00 & $6.38(-8)$  & 3.00 & $2.67(-11)$   & 3.50  &  $1.72(-09)$ & 4.00  \\
$ 256\times 256\times 256$ & $2.04(-12)$  & 3.90 & $7.92(-9)$  & 3.01 & $3.71(-12)$   & 3.61  &  $1.08(-10)$ & 4.00  \\
\hline
\end{tabular}\label{table10}
\end{table}

\subsection{Computational efficiency}

\begin{table}[!tbp]
\tabcolsep = 5pt
\caption{Comparison of the number of iterations, CPU times (in seconds) and errors between the EXCMG method and classical multigrid methods with the Gauss-Seidel smoother. Here CPU$_{W_h}$ denotes
 the computational time for constructing of the initial guess $w_h$.} \centering
\begin{threeparttable}
\begin{tabular}{lccccccccccccc}
    \hline
     &\multirow{2}{*}{$\epsilon$} & \multicolumn{3}{c}{V(1,1)} & & \multicolumn{3}{c}{W(2,1)} & & \multicolumn{4}{c}{EXCMG}   \\
    \cline{3-5} \cline{7-9} \cline{11-14}&  & Iters\tnote{1} &   CPU &   $||u_h-u||_{\infty}$  &  & Iters\tnote{2}  &   CPU &   $||u_h-u||_{\infty}$  & & Iters\tnote{3} &   CPU  &   $||u_h-u||_{\infty}$ & CPU$_{w_h}$ \\
\hline
Exam \ref{exam1}& $10^{-14}$& 16     & 46.1      & $8.61(-13)$  &   &   12        &   47.6    & $8.39(-13)$  & & 8        &   12.9     &  $8.49(-13)$ &  0.6 \\
Exam \ref{exam2}& $10^{-14}$& 16     & 46.4      & $1.91(-12)$  &   &   12        &   47.6    & $1.71(-12)$ & & 9        &   12.6      &  $1.81(-12)$ &  0.6 \\
Exam \ref{gaussian}& $10^{-11}$& 13     & 41.5      & $7.83(-11)$  &   &   9        &   39.4    & $7.22(-11)$ & & 8        &   11.8      &  $2.40(-11)$ &  0.6 \\
Exam \ref{exam4}& $10^{-09}$& 72     & 204.3    & $8.27(-10)$  &    &   47        &   182.9    & $2.76(-10)$ & & 8        &   10.8     &  $1.82(-10)$ &  0.6 \\
Exam \ref{exam5}& $10^{-12}$& 14     & 42.8     & $1.77(-12)$  &    &   10         &   41.3     & $1.75(-12)$ & & 9        &   13.8    &  $1.81(-12)$ &  0.6 \\
Exam \ref{exam3}& $10^{-13}$& 15     & 45.9      & $1.91(-12)$  &   &   11        &   46.1    & $1.91(-12)$ & & 9        &   13.3      &  $2.04(-12)$ &  0.6 \\
\hline
\end{tabular}\label{table11}
\begin{tablenotes}
        \footnotesize
        \item[1] Iters denotes the number of V-cycles required to reach the  error tolerance $\epsilon$.
        \item[2] Iters denotes the number of W-cycles required to reach the  error tolerance $\epsilon$.
        \item[3] Iters denotes the number of CG iterations on the finest grid  for EXCMG method. \end{tablenotes}
          \end{threeparttable}
\end{table}

\begin{table}[!tbp]
\tabcolsep = 5pt
\caption{Comparison of the number of iterations, CPU times (in seconds) and errors between the EXCMG method and classical multigrid methods with the CG smoother. } \centering
\begin{tabular}{lcccccccccccc}
    \hline
     &\multirow{2}{*}{$\epsilon$} & \multicolumn{3}{c}{V(1,1)} & & \multicolumn{3}{c}{W(2,1)} & & \multicolumn{3}{c}{EXCMG}   \\
    \cline{3-5} \cline{7-9} \cline{11-13}&  & Iters &   CPU &   $||u_h-u||_{\infty}$  &  & Iters   &   CPU &   $||u_h-u||_{\infty}$  & & Iters &   CPU  &   $||u_h-u||_{\infty}$ \\
\hline
Exam \ref{exam1}& $10^{-14}$& 15     & 43.9      & $8.38(-13)$  &   &   13        &   49.7    & $8.42(-13)$  & & 8        &   12.9     &  $8.49(-13)$ \\
Exam \ref{exam2}& $10^{-14}$& 15     & 43.3      & $1.76(-12)$  &   &   13        &   49.0    & $1.73(-12)$ & & 9        &   12.6      &  $1.81(-12)$ \\
Exam \ref{gaussian}& $10^{-11}$& 11     & 32.2      & $7.09(-11)$  &   &   10        &   38.0    & $7.22(-11)$ & & 9        &   11.8      &  $2.40(-11)$ \\
Exam \ref{exam4}& $10^{-09}$& 101    & 295.3    & $1.74(-10)$  &    &   30        &   115.8    & $1.71(-10)$ & & 8        &   10.8     &  $1.82(-10)$ \\
Exam \ref{exam5}& $10^{-12}$& 13     & 39.1     & $1.75(-12)$  &    &   11         &   42.7     & $1.75(-12)$ & & 9        &   13.8    &  $1.81(-12)$ \\
Exam \ref{exam3}& $10^{-13}$& 14     & 40.9      & $1.97(-12)$  &   &   11        &   42.4    & $1.95(-12)$ & & 9        &   13.3      &  $2.04(-12)$ \\
\hline
\end{tabular}\label{table12}
\end{table}

In this subsection, we compare the efficiency of the our new EXCMG method with the efficiency of the classical V-cycle and W-cycle multigrid methods for above six examples. Results with Gauss-Seidel smoother are listed in table~\ref{table11} while results with  CG smoother are listed in table~\ref{table12}. In both tables, the number of iterations, computational time, the $L^{\infty}$-norm of the difference between the FD solution $u_h$ and the exact solution $u$ are provided. Moreover, the computational time for constructing of the initial guess $w_h$ (line 6 in algorithm~\ref{alg:oEXCMG}) is also listed in the last column of table~\ref{table11}, this step contains the extrapolation and quartic interpolation as described in section~\ref{3D}.  By comparing the total computational time  of the new EXCMG method with the classical V-cycle and W-cycle multigrid methods for all above six examples as listed in both table~\ref{table11} and~\ref{table12}, one can easily find that the new EXCMG method needs the smallest time for all examples, and this is particularly true when using the unequal meshsizes in different directions, see example \ref{exam4}. Thus, the EXCMG method is much more efficient than the classical V-cycle and W-cycle multigrid methods. And from the last column in table~\ref{table11}, one can find that the computational time for constructing the initial guess $w_h$ described in section~\ref{3D} is $0.6$ seconds for every example, which is very cheap, comparing to the total computational time.

Moreover, one can see from table~\ref{table11} and~\ref{table12}  that there is only a few number of iterations are needed on the finest grid for every example in our EXCMG method, because that the initial guess $w_h$ is already an extremely accurate approximation to the FD solution $u_h$. For example, from the last column of table \ref{table2}, we see that the maximum error of the initial guess on the finest grid for example \ref{exam1} is $4.15\times 10^{-12}$, which implies that the number of significant figures of the approximation exceeds 10.  Indeed, from table \ref{table1}-\ref{table10} we see that the extrapolated value $w_h$ served as an initial guess of the CG solver is a fifth-order approximation to the FD solution $u_h$ in $L^2$-norm, which is one order higher than the convergence order of the fourth-order difference solution $u_h$. Thus, the relative effect of how $w_h$ approximates $u_h$ becomes better when mesh is refined, and the number of iterations is reduced most significantly on the finest grid, see a more detailed discussion in~\cite{Pan2015}.

\begin{figure}
   \centering
   \scalebox{0.6}{\includegraphics{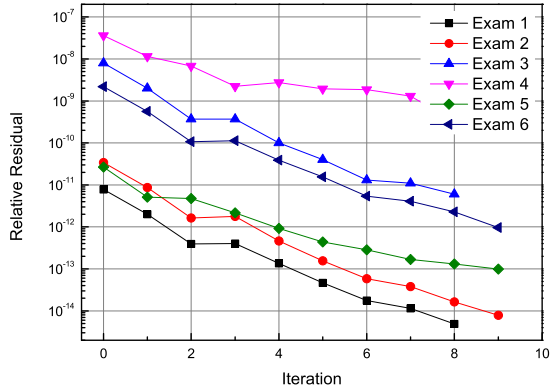}}
   \caption{Relative residual vs. the number of iterations on the finest grid.}\label{RRe}
\end{figure}

Finally, we present the curve of the relative residual on the finest grid versus the number of iterations for the above six examples in Fig.~\ref{RRe}.
As we can see that the initial relative residual on the finest grid for each example is very small. And due to the high oscillations of the initial error as shown in section~\ref{error},
the relative residual decreases by several orders of magnitude  after only a few iterations, and then reaches a number that is less than the required tolerance.

\section{Conclusions}
In this work, we developed a  new extrapolation cascadic multigrid (EXCMG) method  combined with 19-point fourth-order compact difference scheme for solving the 3D Poisson equation on rectangular domains.
The major advantage of the method is to use the Richardson extrapolation and tri-quartic Lagrange interpolation techniques for two numerical solutions on two-level of grids (current and previous grids) to obtain a fifth-order approximation $w_h$ to the fourth-order FD solution $u_h$  as the initial guess of the iterative solution on the next finer grid, which greatly reduces the iteration numbers.  When the exact solution $u$ is sufficiently smooth, a sixth-order extrapolated solution $\tilde{u}_h$ on the fine grid can be obtained by using two fourth-order numerical solutions on two scale grids. Moreover, the gradient of solution $\nabla u_h$ can also be computed easily and efficiently  through solving a series of tridiagonal linear systems resulting from the fourth-order compact FD discretization of the derivatives. Finally, numerical results show that our new  extrapolation cascadic multigrid method is much more efficient comparing to the classical V-cycle and W-cycle multigrid method and it is particularly suitable for solving large scale problems.

The work presented in this paper is an extension of our previous work, which is based on the EXCMG method for the 3D elliptic problem with the linear FE discretization~\cite{Pan2015}. In the near future, we  will extend our method to convection-diffusion equations, Helmholtz equations, biharmonic equations, and other related equations.


\begin{acknowledgements}
Kejia Pan was supported by the National Natural Science Foundation of China (Nos. 41474103 and 41204082), the National
High Technology Research and Development Program of China (No. 2014AA06A602), the Natural Science Foundation of Hunan Province of China (No. 2015JJ3148). Dongdong He was supported by the National Natural Science Foundation of China (No. 11402174), the Program for Young Excellent Talents at Tongji University (No. 2013KJ012)  and the Scientific Research Foundation for the Returned Overseas Chinese Scholars, State Education Ministry. Hongling Hu was supported by the National Natural Science Foundation of China (No. 11301176).
\end{acknowledgements}



\end{document}